\documentclass{eptcs}
\usepackage[utf8]{inputenc}
\usepackage{amsmath,amssymb,amsthm,enumerate,stmaryrd,mathpartir}
\usepackage{hyperref}
\usepackage{graphicx,graphbox}
\usepackage{tikz-cd}
\usepackage{csquotes}
\usepackage[strings]{underscore}

\theoremstyle{definition}
\newtheorem{definition}{Definition}

\newtheorem{proposition}{Proposition}

\newtheorem{corollary}{Corollary}
\newtheorem{lemma}{Lemma}
\newtheorem{remark}{Remark}

\newcommand{\A}{\mathbb{A}}

\newcommand{\C}{\mathbb{C}}
\newcommand{\bh}{\mathbf{H}}
\newcommand{\bv}{\mathbf{V}}

\newcommand{\corner}[1]{{\text{}^\ulcorner_\llcorner\!{#1}\!_\lrcorner^\urcorner}}
\newcommand{\ex}[1]{{#1}^{\circ\bullet}}

\newcommand{\spangraph}{\mathsf{Span}(\mathsf{RGraph})}

\newcommand{\alice}{\texttt{Alice}}
\newcommand{\bob}{\texttt{Bob}}

\makeatletter
\providecommand{\leftsquigarrow}{%
  \mathrel{\mathpalette\reflect@squig\relax}%
}
\newcommand{\reflect@squig}[2]{%
  \reflectbox{$\m@th#1\rightsquigarrow$}%
}
\makeatother

\title{Situated Transition Systems}
\author{Chad Nester \thanks{This research was supported by the ESF funded Estonian IT Academy research measure (project 2014-2020.4.05.19-0001).}\institute{Tallinn University of Technology, Tallinn, Estonia}}

\def\titlerunning{Situated Transition Systems}
\def\authorrunning{Chad Nester}

\begin{document}

\maketitle

\begin{abstract}
We construct a monoidal category of open transition systems that generate material history as transitions unfold, which we call situated transition systems. The material history generated by a composite system is composed of the material history generated by each component. The construction is parameterized by a symmetric strict monoidal category, understood as a resource theory, from which material histories are drawn. We pay special attention to the case in which this category is compact closed. In particular, if we begin with a compact closed category of integers then the resulting situated transition systems can be understood as systems of double-entry bookkeeping accounts.

\end{abstract}

\section{Introduction}

Graphs have been used to model the states and state changes (transitions) of systems for hundreds of years \cite{Euler1741}. Today, graphs can be found everywhere in the scientific literature, and entire fields of study are concerned with specific kinds of graph models. In common practice, to model something as a graph is to treat is as a \emph{closed system} --- that is, the surrounding context is ignored by the model. The closed nature of these models is a failure of compositionality: it prevents us from explaining large systems as the combination of smaller components. This sort of compositionality is all but required if our modelling techniques are to apply to the complex systems we encounter in the world. 

A promising compositional approach is the algebra of transition systems with boundary given by the category $\spangraph$ of spans of reflexive graphs \cite{Kat97}. In this formalism, each transition manifests as an event at the boundaries of a system, and composing systems along a common boundary constrains their behaviour to be consistent with the events observed there. This allows us to consider graph models of \emph{open systems}, and to use these as components in the construction of a larger whole. For example, the authors of \cite{Gia20} have constructed a simplified model of the heart system in the $\spangraph$ setting. 

In an unpublished and --- it seems --- largely unknown paper \cite{Kat98}, the category $\spangraph$ is modified to give a category of systems of partita-doppia (double-entry bookkeeping) accounts. These systems have an account balance, which may change as the result of vaule entering or leaving the system during a transition. The resulting category $\mathsf{Accounts}$ allows us to model a system of partita-doppia accounts in context, as one part of a notional system of all accounts. This is more exciting than may be immediately apparent. From \cite{Kat98}:
\begin{displayquote}
  "The aim of accounting is the measurement of a distributed concurrent system, and it is our contention that it is one of the earliest and most successful mathematical theories of concurrency."
\end{displayquote}

The present work arose from a desire to generalize the category $\mathsf{Accounts}$. In a sense, models in $\spangraph$ (indeed, graph models more generally) are detached from any sort of material reality. The states and state transitions are specified, but the material effect of a given sequence of transitions is left informal, specified as vague intuition. In the category of $\mathsf{Accounts}$, transitions come equipped with a material effect on the partita-doppia ledger associated with that system. The abstract, conceptual world of graphs is thus \emph{situated} in the world of accouting.

Our point of departure is to replace the theory of partita-doppia ledgers with an arbitrary \emph{resource theory} (symmetric strict monoidal category) in the sense of \cite{Coe14}. Augmenting our resource theories with \emph{corners} \cite{Nes21} allows us to assign material history to a transition in a compositional way: material history generated by a composite transition system is the composite of material history generated by its components. We call the resulting notion a \emph{situated transition system}, and we show that for any resource theory $\A$ the $\A$-situated transition systems form a monoidal category.

We show that our formalism specializes to capture its inspiration: if we begin with a compact closed category $\mathbb{Z}$ of integers, the category of $\mathbb{Z}$-situated transition systems is a category of systems of partita doppia accounts in the sense of \cite{Kat98}. Further, we show that for any compact closed category $\A$, the catgory of $\A$-situated transition systems is also compact closed. This generalizes the main theorem of \cite{Kat98}, which is that $\mathsf{Accounts}$ is a compact closed category.

\subsection{Contributions and Related Work}

\emph{Related Work}. We credit the resource-theoretic interpretation of monoidal categories and their string diagrams to \cite{Coe14}. String diagrams for monoidal categories are dealt with rigorously in \cite{Joy91}. The use of ``corners'' in single-object double categories to allow the concurrent decomposition of resource transformations is due to \cite{Nes21}. Double categories first appear in \cite{Ehr63}. Free double categories are considered in \cite{Daw02} and again in \cite{Fio08}. The corner structure we use is in fact the structure of a proarrow equipment. The idea of a proarrow equipment first appears in \cite{Woo82}, albeit in a rather different form. Proarrow equipments have subsequently appeared under many names in formal category theory \cite{Shu08,Gra04}. String diagrams for double categories and proarrow equipments are treated precisely in \cite{Mye16}. The original work on the category of spans of reflexive graphs as a setting for modelling concurrent systems is \cite{Kat97}. Our work is directly inspired by earlier efforts to eqiup such models with accounting information \cite{Kat98}. An excellent mathematical exposition of double-entry bookkeeping is \cite{Ell85}. Compact closed categories were introduced in \cite{Joy96}, along with the compact closed category $\mathbb{Z}$ of integers. More on compact closed categories, and specifically on compact closed categories of integers, can be found in \cite{Abr05}. 

\emph{Contributions}. The main contribution of this paper is the notion of situated transition system, accompanied by the construction of the monoidal category $\mathsf{S}(\A)$ of situated transition systems over an arbitrary monoidal category $\A$ (Propositions \ref{prop:situated-category}, \ref{prop:situated-monoidal}). Other contributions are our investigation into the effect of compact closed structure in $\A$ on $\mathsf{S}(\A)$ (Lemmas \ref{lem:reversal}, \ref{lem:horiz-compact-closed}, \ref{lem:situated-compact-closed}), and the observation that $\mathsf{S}(\mathbb{Z})$ captures the systems of partita-doppia accounts of \cite{Kat98} (Corollary \ref{cor:partita-doppia-corollary}). To our knowledge the compact closed perspective on double-entry bookkeeping is also novel, and so may be viewed as a modest contribution.

\section{Preliminaries}

\subsection{Monoidal Categories as Resource Theories}\label{sec:resourcetheories}

Symmetric strict monoidal categories can be understood as theories of resource transformation \cite{Coe14}. Objects are interpreted as collections of resources, with $A \otimes B$ the collection consisting of both $A$ and $B$, and $I$ the empty collection. Arrows $f : A \to B$ are understood as ways to transform the resources of $A$ into those of $B$, or equivalently as parts of a larger \emph{material history} involving those resources. We call symmetric strict monoidal categories \emph{resource theories} when we have this sort of interpretation in mind.

For example, let $\mathfrak{B}$ be the free symmetric strict monoidal category generated by:
\begin{mathpar}
  \{ \mathsf{bread},\mathsf{dough},\mathsf{flour},\mathsf{oven}\}\\
  
  \mathsf{knead} : \mathsf{flour} \to \mathsf{dough}

  \mathsf{bake} : \mathsf{dough} \otimes \mathsf{oven} \to \mathsf{bread} \otimes \mathsf{oven}

  \mathsf{eat} : \mathsf{bread} \to I
\end{mathpar}
subject to no equations. $\mathfrak{B}$ can be understood as a resource theory of bread. The arrow $\mathsf{knead}$ represents the process of making $\mathsf{dough}$ from $\mathsf{flour}$, $\mathsf{bake}$ represents baking $\mathsf{dough}$ in an $\mathsf{oven}$ to obtain $\mathsf{bread}$ (and an $\mathsf{oven}$), and $\mathsf{eat}$ represents the consumption of $\mathsf{bread}$. 

The structure of symmetric strict monoidal categories provides natural algebraic scaffolding for composite transformations, with the associated string diagrams acting as a convenient syntax for expressing material histories. For example in the following string diagram over $\mathfrak{B}$ we see two units of $\mathsf{dough}$ made into loaves of $\mathsf{bread}$ by baking one after the other in an $\mathsf{oven}$. 
\[
\includegraphics[height=3cm]{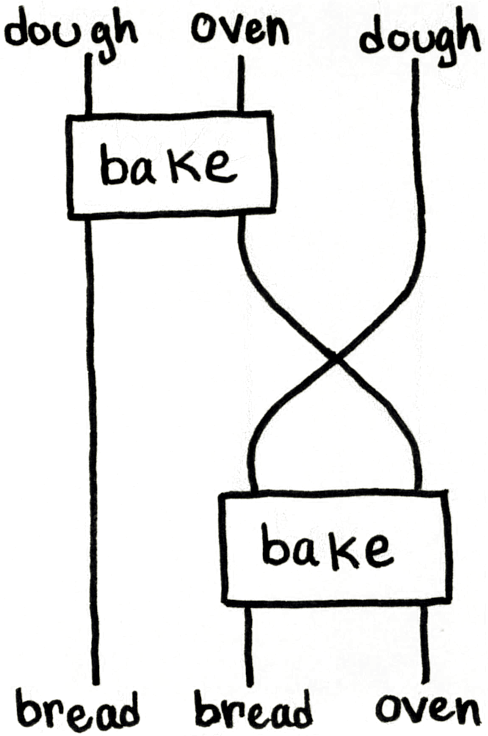}
\]
Notice how the topology of the diagram captures the logical flow of resources.

Given a parallel pair $f,g : A \to B$ of material histories in some resource theory $\A$, we understand equality of $f$ and $g$ to mean that both have the same effect on the resources involved. For example, suppose we add a generating morphism $\mathsf{sift} : \mathsf{flour} \to \mathsf{flour}$ to our resource theory $\mathfrak{B}$, subject to the equation $\mathsf{sift} \circ \mathsf{sift} = \mathsf{sift}$. Call the resulting resource theory $\mathfrak{B}_\mathsf{sift}$. In this new theory the material histories $\mathsf{sift}$ and $\mathsf{sift} \circ \mathsf{sift}$ express different sequences of events, with the $\mathsf{flour}$ being sifted once in the former, but twice in the latter. They are made equal by our new equation, which means that in $\mathfrak{B}_\mathsf{sift}$, sifting $\mathsf{flour}$ twice has the same effect as sifting it once. Contrast this to $1_\mathsf{flour}$ and $\mathsf{sift} : \mathsf{flour} \to \mathsf{flour}$. Identity morphisms have no effect on the resources involved, so intuitively these two material histories should not denote equal morphisms of $\mathfrak{B}_\mathsf{sift}$, and indeed they do not. We adopt this understanding of equality as a general principle in our design and understanding of resource theories. 

\subsection{Cornering and Concurrent Transformations}

The resource theoretic interpretation of symmetric strict monoidal categories can be extended to allow the decomposition of material histories into their concurrent components \cite{Nes21}. Specifically, we augment the string diagrams for a given resource theory $\A$ with \emph{corners} for each object $A$ of $\A$:
\[
\includegraphics[height=1.7cm,align=c]{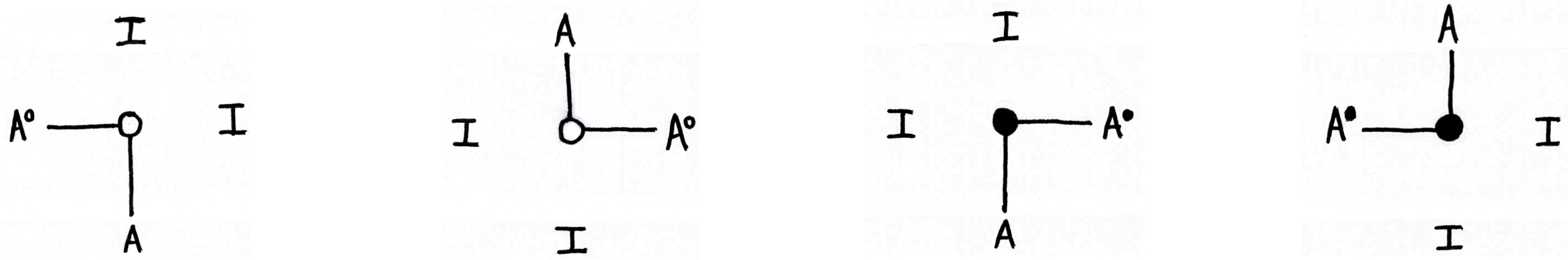}
\]
Corners allow us to express resources flowing into and out of a system. $A^\circ$ denotes an instance of $A$ flowing from left to right, and $A^\bullet$ denotes an instance of $A$ flowing from right to left. Our corners must satisfy the \emph{yanking identities}, which ensure that this movement has no effect on the resources themselves:
\[
\includegraphics[height=1cm,align=c]{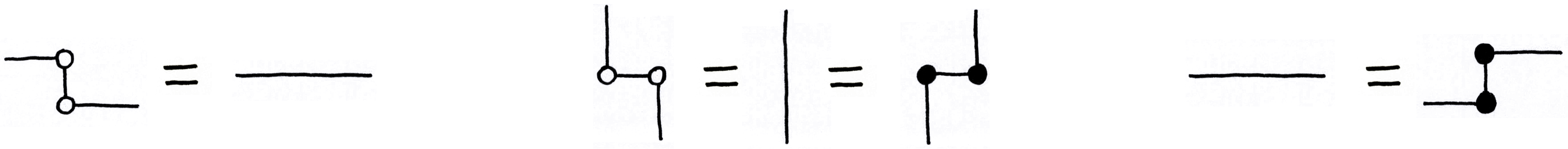}
\]
For example, adding corners to our resource theory $\mathfrak{B}$ allows the following decomposition of the baking process. The transformation below on the left begins with no resources, then $\mathsf{flour}$ enters along the right boundary and is $\mathsf{knead}$ed into $\mathsf{dough}$, which leaves along the right boundary. The transformation below in the middle begins with an $\mathsf{oven}$, then $\mathsf{flour}$ passes through from right to left, $\mathsf{dough}$ is received along the left boundary and is $\mathsf{bake}$d, and the resulting $\mathsf{bread}$ leaves along the right boundary, with the $\mathsf{oven}$ staying put. Finally, the transformation below on the right begins with $\mathsf{flour}$, which leaves the system along the left boundary, after which $\mathsf{bread}$ enters from the left, and is $\mathsf{eat}$en. 
\begin{mathpar}
\includegraphics[height=1.7cm,align=c]{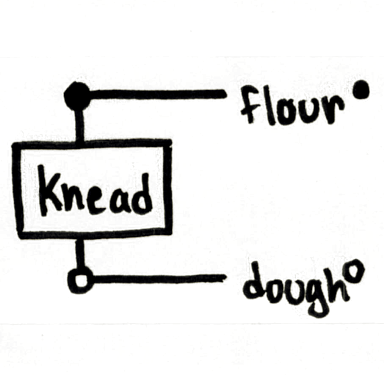}

\includegraphics[height=1.7cm,align=c]{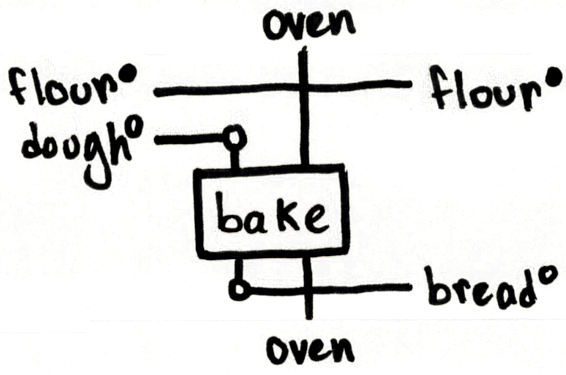}
    
\includegraphics[height=1.7cm,align=c]{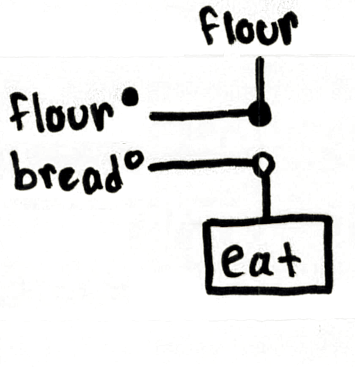}
\end{mathpar}
These transformations may be composed horizontally to obtain a single transformation of resources:
\[
\includegraphics[height=3cm,align=c]{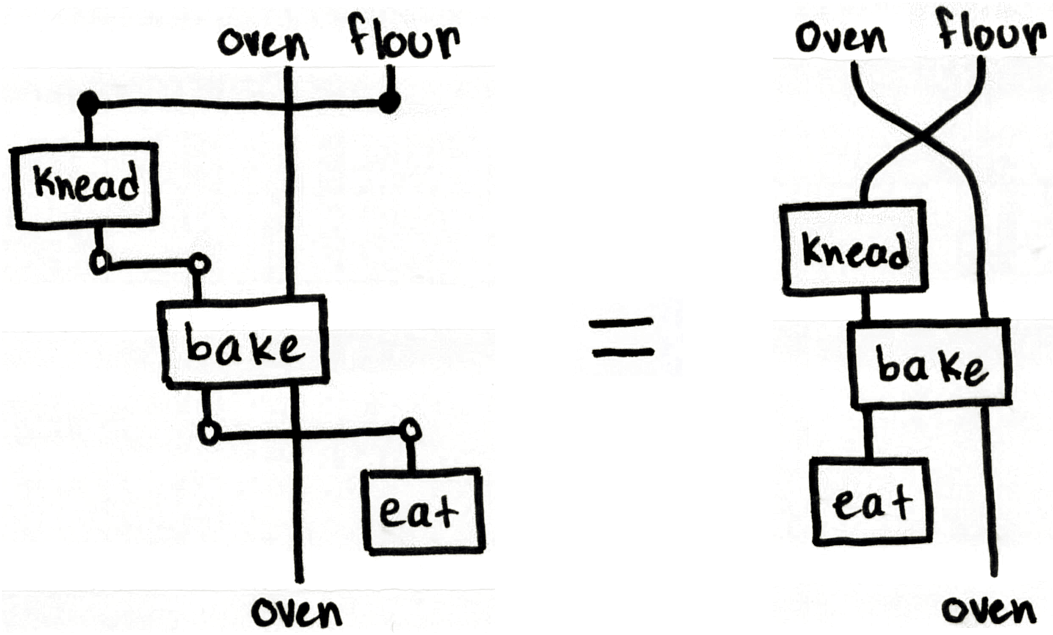}
\]
Formally, these augmented string diagrams denote cells of a single-object double category $\corner{\A}$ which we call the \emph{free cornering of $\A$}. This double category has one object, so in particular the horizontal and vertical edge categories are necessarily monoids (single-object categories). The horizontal edge monoid $(\A_0,\otimes,I)$ is given by the monoidal structure on the objects of $\A$. The vertical edge monoid $\ex{\A} = (\A_0 \times \{\circ, \bullet\})^*$ is the free monoid of polarized objects of $\A$, written as in $A^\circ$ and $A^\bullet$. Elements of $\ex{\A}$ are sequences of polarized objects of $\A$, which we understand as \emph{$\A$-valued exchanges}. The monoid operation is given by concatenation (denoted by $\otimes$) and the empty sequence (denoted by $I$) is the unit of the monoid. Each exchange $X_1 \otimes \cdots \otimes X_n \in \ex{\A}$ involves a left participant and a right participant giving each other resources in sequence, with $A^\circ$ indicating that the left participant should give the right participant an instance of $A$, and $A^\bullet$ indicating that the right participant should give the left participant an instance of $A$. For example if $\alice$ is the left participant and $\bob$ is the right participant, then we can picture the exchange $A^\circ \otimes B^\bullet \otimes C^\bullet \in \ex{\A}$ as 
\[
\alice \rightsquigarrow
\includegraphics[height=1cm,align=c]{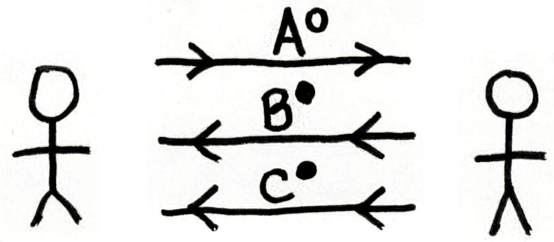}
\leftsquigarrow \bob
\]
These exchanges happen in order. The exchange pictured above demands that first $\alice$ gives $\bob$ an instance of $A$, then $\bob$ gives $\alice$ an instance of $B$, and then finally $\bob$ gives $\alice$ an instance of $C$.

The generating cells of $\corner{\A}$ are the corners discussed above, subject to the yanking equations, together with cells $\corner{f}$ for each arrow $f: A \to B$ of $\A$, subject to the following equations: 
\[
  \includegraphics[height=1.7cm,align=c]{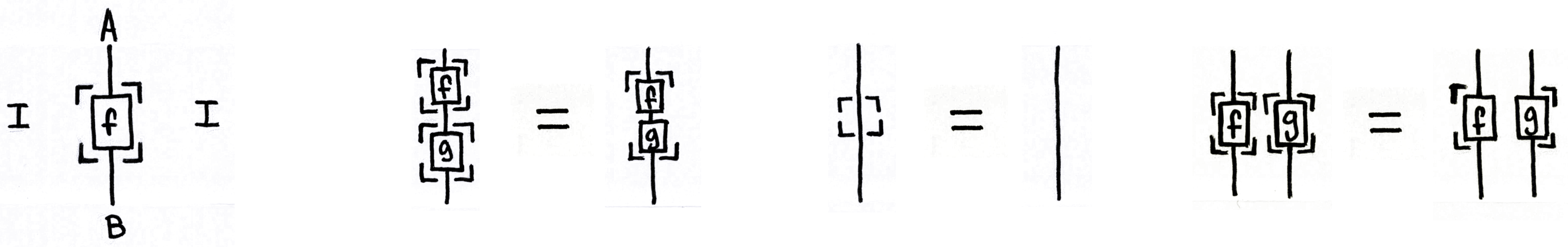}
\]
Now $\corner{\A}$ is the free double category generated by this data, with arbitrary cells of $\corner{\A}$ being obtained by vertical and horizontal composition of the generators, subject to the equations of a double category (see \cite{Fio08,Daw02} for more on free double categories).

The double category $\corner{\A}$ is more thoroughly investigated in \cite{Nes21}. For our purposes we need only mention that $\corner{\A}$ always contains \emph{crossing cells}, pictured below on the left for an arbitrary $B \in \A_0$ and $X \in \ex{\A}$. These crossing cells make $\corner{\A}$ into a monoidal double category in the sense of \cite{Shu10}, with the tensor product of cells given given below on the right.
\begin{mathpar}
\includegraphics[height=1.7cm]{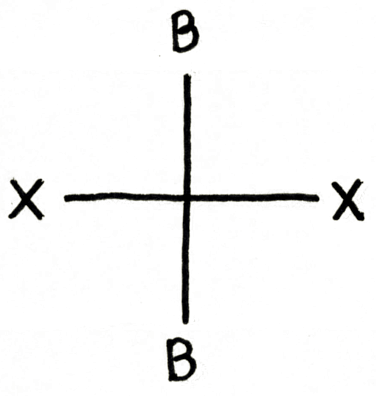}

\includegraphics[height=1.7cm]{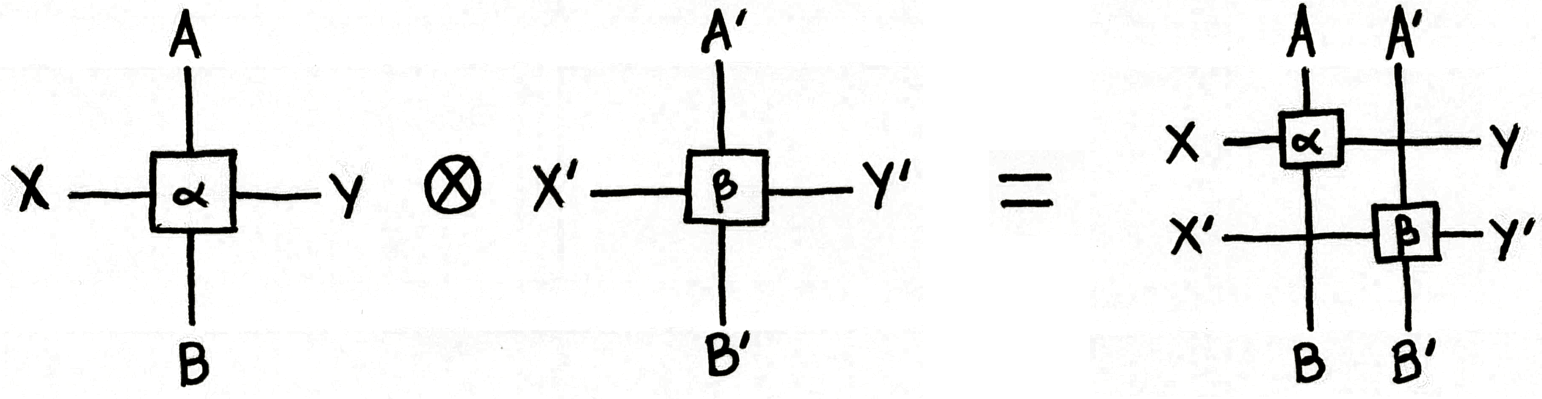}
\end{mathpar}

This is all the resource-theoretic machinery we will need to give a compositional account of the material histories generated by our transition systems. We turn now to the transition systems themselves.

\subsection{The Algebra of Transition Systems with Boundary} 

For our purposes a \emph{transition system} $R$ consists of a collection of \emph{states}, $R_0$, and a collection of \emph{transitions} $t : A \to B \in R_1$ where $A,B \in R_0$. We ask further that for each $A \in R_0$ there is a \emph{trivial transition} $\varepsilon_A : A \to A \in R_1$. In other words, a transition system is precisely a \emph{reflexive graph} (states are vertices, transitions are edges). A \emph{morphism} $F : R \to S$ of transition systems is a morphism of reflexive graphs: It consists of a mapping of vertices $F_0 : R_0 \to S_0$ together with a mapping of edges $F_1 : R_1 \to S_1$ and must preserve the source and target of edges in the sense that if $t : A \to B$ then $F_1(t) : F_0(A) \to F_0(B)$. Further, it must preserve the trivial edges in the sense that $F_1(\varepsilon_A) = \varepsilon_{F_0(A)}$. Reflexive graphs and reflexive graph morphisms form a cartesian category $\mathsf{RGraph}$, which will play a supporting role in our development.

The algebra of transition systems with boundary is captured by the category $\spangraph$ of spans in $\mathsf{RGraph}$ \cite{Kat97}. If $U$ and $V$ are reflexive graphs, then a morphism $R : U \to V$ of $\spangraph$ consists of another reflexive graph $R$ (the \emph{apex}) with morphisms $\delta_0 : R \to U$ and $\delta_1 : R \to V$ of $\mathsf{RGraph}$ (the \emph{legs}). 
We understand this as a transition system $R$ with \emph{boundaries} $U$ and $V$. Every transition $t : A \to B$ of $R$ corresponds to an \emph{event} at each boundary --- $\delta_0(t)$ at $U$ and $\delta_1(t)$ at $V$. Span composition is given by pullback: If $R : U \to V$ and $S : V \to W$ in $\spangraph$, a transition of $S \circ R : U \to W$ consists of a pair of transitions $(t,t') \in R_1 \times S_1$ which correspond to the same event $\delta_1(t) = \delta_0(t')$ at the shared boundary $V$. In the composite each of the components constrains the behaviour of the other. We consider spans modulo the equivalence relation generated by span isomorphism.

For example, let $M$ be the reflexive graph with a single vertex and two nontrivial edges $\mathsf{up}$ and $\mathsf{down}$, pictured below on the left. The diagram below on the right indicates a morphism $\mathsf{Gear} : M \to M$ of $\spangraph$. The apex has a single vertex and two nontrivial edges $\mathsf{cw}$ and $\mathsf{ccw}$, and the legs of the span are indicated by the colouring. The idea is that our gear can rotate clockwise ($\mathsf{cw}$), in which case the teeth along the left and right boundary move $\mathsf{up}$ and $\mathsf{down}$ respectively, or may rotate counterclockwise ($\mathsf{ccw}$), with the boundary teeth moving in the opposite directions. We omit the trivial edges from our diagrams but nonetheless consider them to be present, so our gear system can also do nothing via $\varepsilon$. 
\begin{mathpar}
{M = \includegraphics[height=1.2cm,align=c]{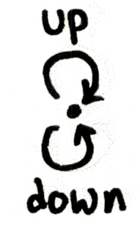}}

{\mathsf{Gear} : M \to M \hspace{0.5cm} = \hspace{0.5cm}
\begin{array}{c|c|c}
  \includegraphics[height=1.2cm,align=c]{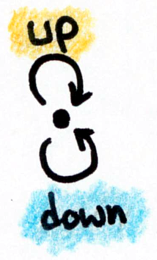}
  &
  \includegraphics[height=1.2cm,align=c]{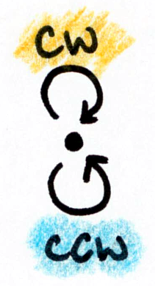}
  &
  \includegraphics[height=1.2cm,align=c]{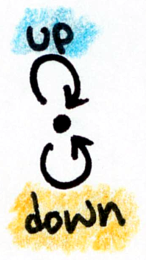}
\end{array}}
\end{mathpar}

Now the composite system $\mathsf{Gear} \circ \mathsf{Gear}$ represents two interlocking gears. The teeth interlock at the shared boundary, where they must move in unison. Our notion of composition captures this formally: the apex of our composite span has a single vertex and two nontrivial edges, one in which the gear on the left rotates clockwise and the gear on the right rotates counterclockwise, and one representing the opposite situation. The case where both gears rotate in the same direction is not present as it would be inconsistent along the shared boundary. In fact $\mathsf{Gear} \circ \mathsf{Gear} = 1_M$, reflecting a similar property of physical gears.

\[
\mathsf{Gear} \circ \mathsf{Gear} : M \to M \hspace{0.5cm} = \hspace{0.5cm}
\begin{array}{c|c|c}
  \includegraphics[height=1.2cm,align=c]{figs/gear-left-boundary.png}
  &
  \includegraphics[height=1.2cm,align=c]{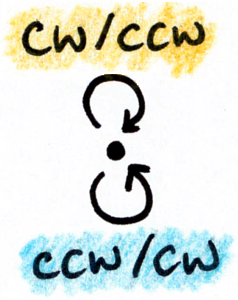}
  &
  \includegraphics[height=1.2cm,align=c]{figs/gear-left-boundary.png}
\end{array}
\]

$\spangraph$ is a symmetric monoidal category. The tensor product is defined on objects by $U \otimes V = U \times V$, and the unit $1$ is the graph with a single vertex and no nontrivial edges. On arrows $R : U \to V$ and $S : U' \to V'$ the tensor product $R \otimes S : U \otimes U' \to V \otimes V'$ has apex $R \times S$ with left and right leg given by the product of the left and right legs of $R$ and $S$, respectively. A transition in the tensor product of two systems is simply a transition from each component. Intuitively, the components function independently of each other. Further, notice that the component systems may function asynchronously via the $\varepsilon$ transitions: If $t \in R_1$ and $t' \in S_1$ then $(t,t'), (t,\varepsilon), (\varepsilon,t')$, and $(\varepsilon,\varepsilon)$ are all transitions of $R \otimes S$. 

There is also a lot of other structure in $\spangraph$. Relevant to our purposes here is the fact that $\spangraph$ is compact closed. The dual of $X$ is given by $X$ itself, and the unit and counit are defined in terms of the finite product structure on $\mathsf{RGraph}$: $\eta_X : 1 \to X \otimes X$ is given by the span with apex $X$, left leg $!_{X} : X \to 1$, and right leg $\Delta_X : X \to X \times X$, with $\varepsilon_X : X \otimes X \to 1$ constructed similarly.

We conclude our discussion of $\spangraph$ with a bread-themed example. Define objects $U,V$ of $\spangraph$ as follows --- again omitting the trivial edges from our diagrams:
\begin{mathpar}
U = \includegraphics[height=0.5cm,align=c]{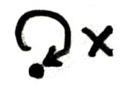}

V = \includegraphics[height=0.5cm,align=c]{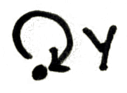}
\end{mathpar}
We understand the event $x \in U_1$ to indicate that the system on the right is obtaining ingredients for baking from the system on the left, and the $y \in V_1$ indicates that the system on the left is selling bread to the system on the right.

Let $\mathsf{Baker}$ be the morphism of $\spangraph$ pictured below on the left. The apex has two vertices, one in which the system is $\mathsf{open}$ for business, and another in which it is $\mathsf{closed}$. There are edges allowing the system to transition from being $\mathsf{open}$ to being $\mathsf{closed}$, and vice versa. When it is $\mathsf{open}$, the system may $\mathsf{bake}$ and $\mathsf{sell}$ bread. The legs of the span are indicated by the colouring: The $\mathsf{bake}$ transition corresponds to the event $x$ at the left boundary, and the transition $\mathsf{sell}$ corresponds to the event $y$ at the right boundary. An absence of colour indicates the trivial event $\varepsilon$, so for example the transition $\mathsf{open}$ corresponds to the trivial event at both boundaries, and $\mathsf{bake}$ corresponds to the trivial event at the right boundary. 
\begin{mathpar}
{\mathsf{Baker} : U \to V \hspace{0.1cm} =
\begin{array}{c|c|c}
  \includegraphics[height=0.5cm,align=c]{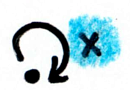}
  &
  \includegraphics[height=1.7cm,align=c]{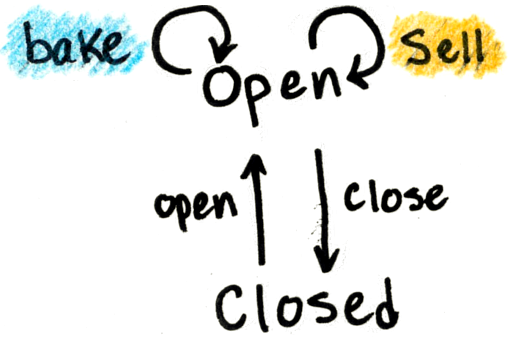}
  &
  \includegraphics[height=0.5cm,align=c]{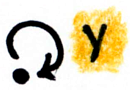}
\end{array}}

{\mathsf{Eater} : V \to 1 \hspace{0.1cm} = 
\begin{array}{c|c}
  \includegraphics[height=0.5cm,align=c]{figs/y-loop-yellow.png}
  &
  \includegraphics[height=1.8cm,align=c]{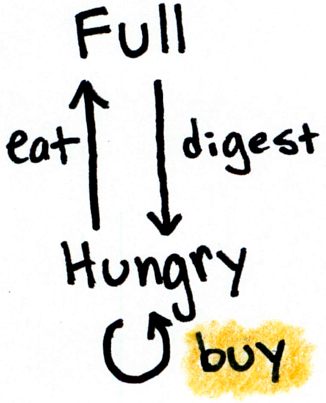}
\end{array}}
\end{mathpar}
Let $\mathsf{Eater}$ be the morphism of $\spangraph$ pictured above on the right. The apex has two vertices, one in which the system is $\mathsf{hungry}$, and another in which it is $\mathsf{full}$. If $\mathsf{hungry}$, the system may $\mathsf{eat}$ to become $\mathsf{full}$, and if $\mathsf{full}$ may $\mathsf{digest}$ to become $\mathsf{hungry}$. Finally, when it is $\mathsf{hungry}$ the system may $\mathsf{buy}$ food. The legs are again indicated by the colouring, with the right leg omitted entirely since in this case there is nothing to indicate. The transition $\mathsf{buy}$ corresponds to event $\mathsf{y}$ at the left boundary, and that is all. 

Now, composing our two systems along their shared boundary $V$ yields:
\begin{mathpar}
{\mathsf{Eater} \circ \mathsf{Baker} : U \to 1 \hspace{0.1cm} =
\begin{array}{c|c}
  \includegraphics[height=0.5cm,align=c]{figs/x-loop-blue.png}
  &
  \includegraphics[height=2.4cm,align=c]{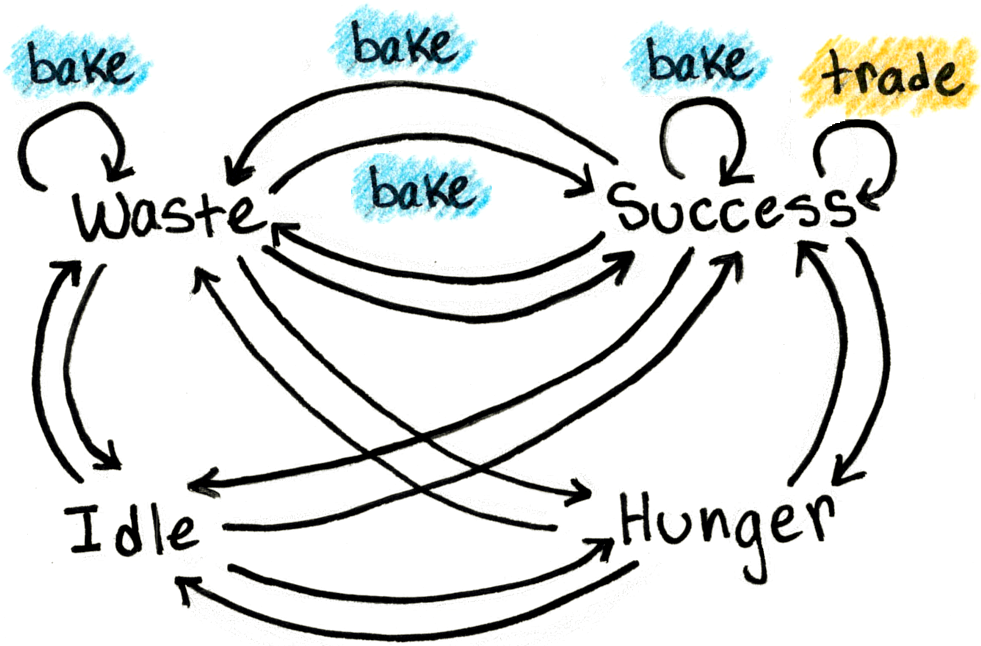}
\end{array}}

{\includegraphics[height=1.3cm,align=c]{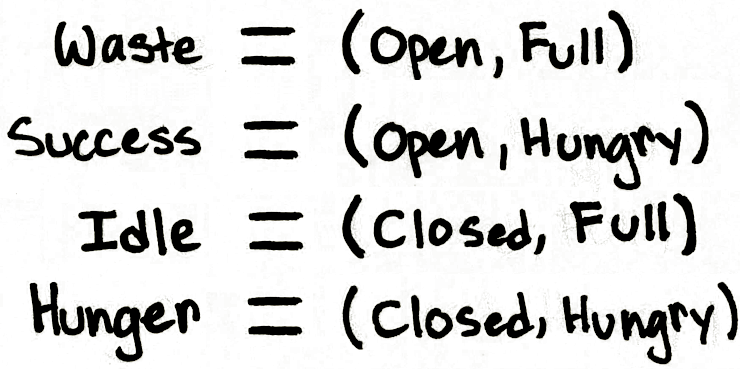}}
\end{mathpar}
The unlabelled transitions arise from combinations of $\mathsf{open}$, $\mathsf{close}$, $\mathsf{eat}$, and $\mathsf{digest}$ --- those transitions corresponding to the trivial event at the boundaries. The $\mathsf{bake}$ transitions are those in which the $\mathsf{Baker}$ system $\mathsf{bake}$s, and the $\mathsf{trade}$ transition corresponds to the $\mathsf{Baker}$ subsystem $\mathsf{sell}$ing bread and the $\mathsf{Eater}$ subsystem $\mathsf{buy}$ing it --- activities which must be synchronised in the composite system. The legs of the span are indicated by the colouring, and we see that every $\mathsf{bake}$ transition involves the event $x$ along the left boundary. The transition $\mathsf{trade}$ is coloured yellow to draw attention to the fact that it is the coincidence of the two yellow transitions in the component systems, and it has trivial boundary events.

\section{Situated Transition Systems}

Given a resource theory $\A$, in this section we show how transition systems with boundary can be equipped to generate $\A$-valued material histories as transitions occur. The double category $\corner{\A}$ of concurrent transformations plays an essential role, allowing us to combine the histories generated by component spans into the history generated by their composite through horizontal composition in $\corner{\A}$. 

We begin by situating the boundaries of our transition systems. In $\spangraph$ the possible events (edges) along a boundary (reflexive graph) serve to synchronise and constrain the behaviour of the larger system. From the material point of view, the relevant part of a boundary event is whether or not any resources leave or enter the system, and if so which ones. This information is captured by the monoid $\ex{\A}$ of $\A$-valued exchanges, which is equivalently a reflexive graph with a single vertex where the unit $I$ of the monoid is the trivial edge. 
\begin{definition}
  Let $\A$ be a resource theory. Then an \emph{$\A$-situated boundary} $(U,\phi_U)$ consists of a reflexive graph $U$ together with a reflexive graph homomorphism $\phi_U : U \to \ex{\A}$. Call $\phi_U$ the \emph{situation of $U$ in $\A$}. 
\end{definition}
We understand $\phi_U(x)$ to describe the resources that cross the boundary as part of the event $x$, and thus constitute its material effect. We will depict $\A$-situated boundaries as graphs with edge labels drawn from $\ex{\A}$, defining the situation of the boundary in $\A$. Since $\ex{\A}$ has only one vertex, we do not need to label the vertices. Edges with no label are understood as having label $I$, and we continue to omit the trivial edges from our depictions. For $X \in \ex{\A}$ we adopt the convention of writing $X$ for the $\A$-situated boundary with a single vertex and a single nontrivial edge, which is mapped to $X$ by the situation. For example the $\mathfrak{B}$-situated boundary $\mathsf{flour}^\circ$ is depicted below on the left. The boundary with two vertices and two nontrivial edges --- one from each vertex to the other --- which are both mapped to $I$ by the situation is depicted below on the right. 
\begin{mathpar}
  \includegraphics[height=.6cm,align=c]{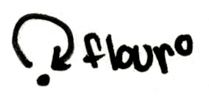}

  \includegraphics[height=.6cm,align=c]{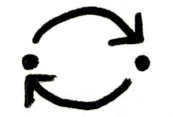}
\end{mathpar}

Now to situate entire transition systems we associate each transition with a cell of $\corner{\A}$ describing the corresponding material effect. The left and right boundaries of this cell must match the labels in $\ex{\A}$ of the left and right boundary events, respectively, so that any material exchanges entailed by those events are present in the material history of the transition. In order to make this precise we view $\corner{\A}$ as a span of reflexive graphs. Specifically, define $\left< \A \right>$ to be the reflexive graph with vertex set $\A_0$ in which an edge $\alpha : A \to B$ is a cell $\alpha$ of $\corner{\A}$ with top boundary $A$ and bottom boundary $B$. Then there is a span
\[ \begin{tikzcd}
\ex{\A} & \left< \A \right> \ar[r,"\delta_1"] \ar[l,"\delta_0"'] & \ex{\A}
\end{tikzcd} \]
where $\delta_0(\alpha)$ and $\delta_1(\alpha)$ are the left and right boundary of $\alpha$, respectively. The trivial edges of $\left< \A \right>$ are given by the vertical identity cells. Situated transition systems are now defined as follows.

\begin{definition}
  Let $\A$ be a resource theory, and let $(U,\phi_U)$ and $(V,\phi_V)$ be $\A$-situated boundaries. Then an \emph{$\A$-situated transition system} $(R,\phi_R) : (U,\phi_U) \to (V,\phi_V)$ consists of a morphism $U \leftarrow R \to V$ of $\spangraph$ together with a reflexive graph homomorphism $\phi_R : R \to \left< \A \right>$ that we call the \emph{situation of $R$ in $\A$}. We require $\phi_R$ to be \emph{coherent} with respect to $\phi_U$ and $\phi_V$ in the sense that the following diagram of reflexive graph homomorphisms commutes:
  \[
  \begin{tikzcd}
    U \ar[d,"\phi_U"'] & R \ar[d,"\phi_R"] \ar[l] \ar[r] & V \ar[d,"\phi_V"] \\
    \ex{\A} & \left< \A \right> \ar[l,"\delta_0"] \ar[r,"\delta_1"'] & \ex{\A}
  \end{tikzcd}
  \]
\end{definition}
We understand $\phi_R$ as assigning a collection of resources to each state of $R$, and assigning to each transition of $R$ a concurrent transformation of resources whose left and right boundary coincide with the material effect of the left and right boundary events. We depict situated transition systems by giving the underlying span of reflexive graphs as before, with the legs indicated by the colouring. We indicate the action of $\phi_R$ by labelling the vertices (resp.\ edges) of the apex with the object of $\A$ (resp.\ cell of $\corner{\A}$) that $\phi_R$ maps them to. For example we can refine our earlier bread-themed example to be $\mathfrak{B}$-situated, with the new $\mathsf{Baker}$ system given by:

\[
\mathsf{Baker} : \mathsf{flour}^\circ \to \mathsf{bread}^\circ \hspace{0.5cm} = \hspace{0.5cm}
\begin{array}{c|c|c}
  \includegraphics[height=0.5cm,align=c]{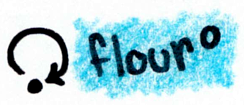}
  &
  \includegraphics[height=1.7cm,align=c]{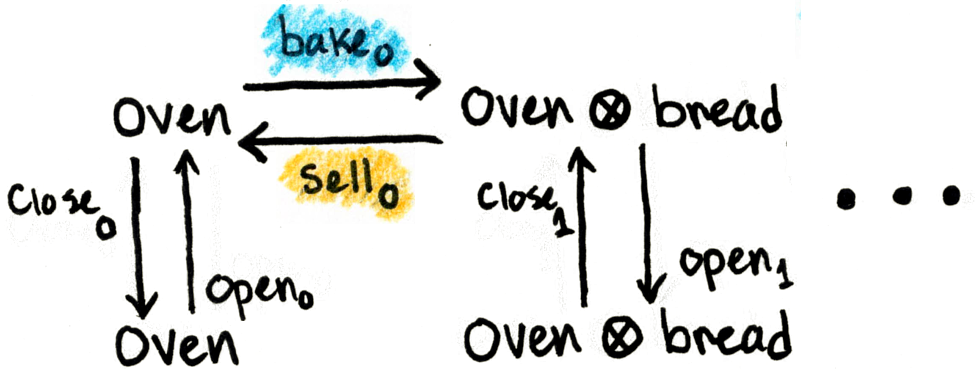}
  &
  \includegraphics[height=0.5cm,align=c]{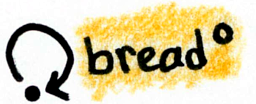}
\end{array}
\]
\noindent where the edge labels are the following cells of $\corner{\mathfrak{B}}$:
\begin{mathpar}
  \mathsf{bake}_n = \includegraphics[height=2.7cm,align=c]{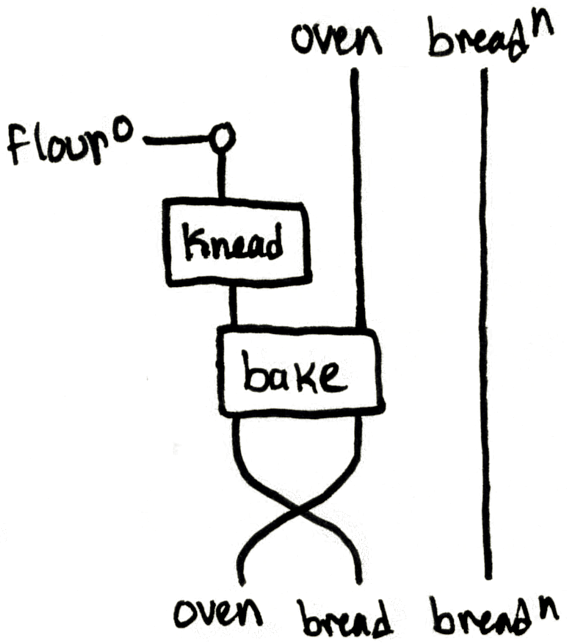}

  \mathsf{sell}_n = \includegraphics[height=1.7cm,align=c]{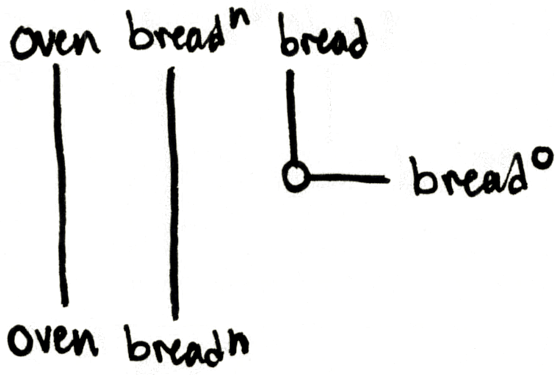}

  \mathsf{close}_n = \mathsf{open}_n = \includegraphics[height=1.7cm,align=c]{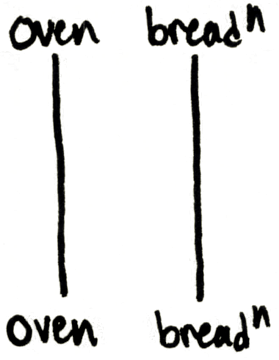}
  
\end{mathpar}
The left boundary is given by the graph with a single vertex and one nontrivial edge, which is mapped to $\mathsf{flour}^\circ$ by the situation, indicating that $\mathsf{flour}$ enters the system as part of that event. The right boundary is similar, with the single nontrivial edge mapped to $\mathsf{bread}^\circ$ by the situation, indicating that $\mathsf{bread}$ leaves the system. The apex has two vertices for each $n \in \mathbb{N}$ which indicate whether the system is open for business or not, and that it currently has $n$ units of $\mathsf{bread}$ in stock. The two states in which the sytem has $n$ units of bread are mapped to $\mathsf{oven} \otimes \mathsf{bread}^n$ by the situation. The edges are similarly indexed: the system may open and close while retaining its stores of bread via $\mathsf{open}_n$ and $\mathsf{close}_n$. When open the system may bake bread via $\mathsf{bake}_n$, in which case we see that $\mathsf{flour}$ enters the system from the left, and may also sell any bread it has via $\mathsf{sell}_n$, in which case $\mathsf{bread}$ leaves from the right.

We continue by defining a $\mathfrak{B}$-situated $\mathsf{Eater}$ as follows:
\[
\mathsf{Eater} : \mathsf{bread}^\circ \to I \hspace{0.5cm} = \hspace{0.5cm}
\begin{array}{c|c}
  \includegraphics[height=0.5cm,align=c]{figs/bread-loop-yellow.png}
  &
  \includegraphics[height=1.5cm,align=c]{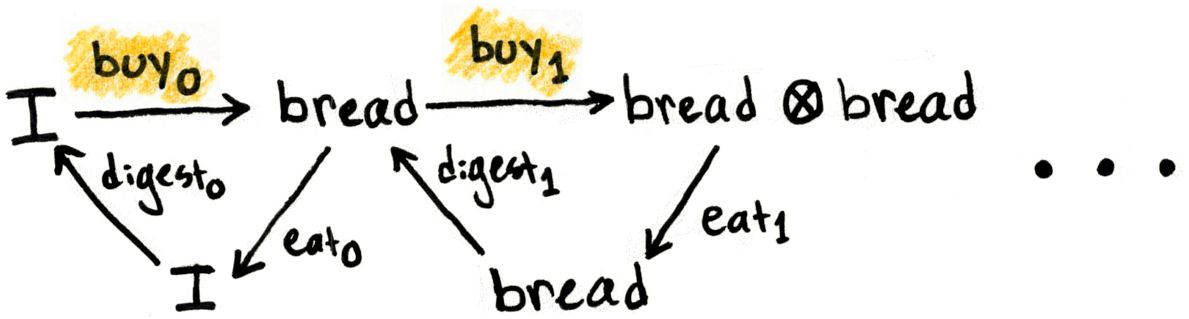}
\end{array}
\]
where the edge labels are the following cells of $\corner{\mathfrak{B}}$:
\begin{mathpar}
  \mathsf{buy}_n = \includegraphics[height=1.7cm,align=c]{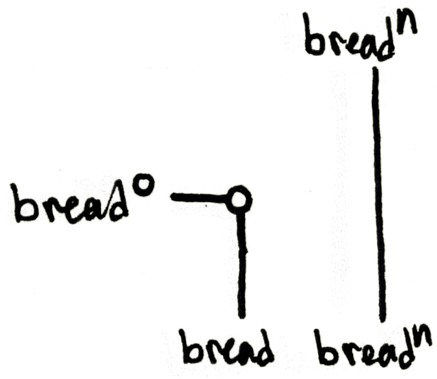}

  \mathsf{eat}_n = \includegraphics[height=1.7cm,align=c]{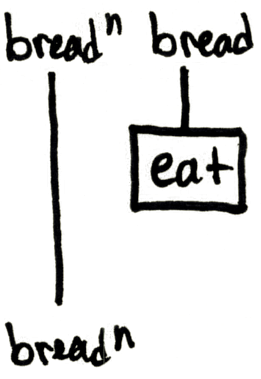}

  \mathsf{digest}_n = \includegraphics[height=1.7cm,align=c]{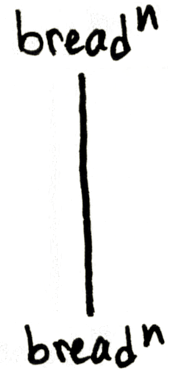}
\end{mathpar}
There are two states for each $n \in \mathbb{N}$ in which the system is hungry, and one in which it is full. In the $n$th iteration of each of these states, the system posesses $n$ units of $\mathsf{bread}$. If in a hungry state and posessing at least one $\mathsf{bread}$, the $\mathsf{eat}_n$ transitions allow it to eat and enter a full state. From a full state the $\mathsf{digest}_n$ transitions allow the system to become hungry, leaving the amount of bread unchanged, and finally if the system is hungry then the $\mathsf{buy}_n$ transitions allow it to acquire more $\mathsf{bread}$ along the left boundary, with the legs of the span indicating that when this happens $\mathsf{bread}$ must enter the system along the left boundary. 

To compose $\A$-situated transition systems $(R,\phi_R) : (U,\phi_U) \to (V,\phi_V)$ and $(S,\phi_S) : (V,\phi_V) \to (W,\phi_W)$ we compose the underlying spans by pullback as in $\spangraph$, and define the composite situation $\phi_{S \circ R} : S \circ R \to \left< \A \right> $ by horizontal composition: $\phi_{S \circ R}(t,t') = \phi_R(t) \mid \phi_S(t')$. This is well-defined because the situations are coherent. In particular this means that $\delta_1 \circ \phi_R = \delta_0 \circ \phi_S$, which says precisely that the right boundary of $\phi_R(t)$ is the left boundary of $\phi_S(t')$ for edges $(t,t')$ of $S \circ R$. Composition of situated transition systems is associative because composition in $\spangraph$ and horizontal composition in $\corner{\A}$ are both associative. Notice also that paths in a situated transition system have vertically composable material effects, with the composite giving the effect of the entire sequence of transitions.

Continuing our example, we may compose our $\mathfrak{B}$-situated $\mathsf{Eater}$ and $\mathsf{Baker}$ transition systems to obtain $\mathsf{Eater} \circ \mathsf{Baker} : \mathsf{flour}^\circ \to I$. This transition system has four vertices for each pair $n,m$ of natural numbers, being those states in which the $\mathsf{Baker}$ has $n$ $\mathsf{bread}$ and the $\mathsf{Eater}$ has $m$ $\mathsf{bread}$. The transitions of this new system are mostly pairs of transitions of the components, the exception being that when the $\mathsf{Baker}$ $\mathsf{sell}$s the $\mathsf{Eater}$ must $\mathsf{buy}$ due to the fact that these transitions are assigned to the same event along the shared boundary $\mathsf{bread}^\circ$. Now, suppose that in our composite system the $\mathsf{Baker}$ begins with one $\mathsf{bread}$ and that the $\mathsf{Eater}$ begins with none. Suppose further that events unfold as follows: First, the $\mathsf{Baker}$ sells its bread to the $\mathsf{Eater}$, which promptly eats it. Then, the $\mathsf{Baker}$ bakes more bread, and finally sells the new bread to the $\mathsf{Eater}$. This sequence of transitions corresponds to the following material history: below on the left we see the history generated by the $\mathsf{Baker}$, below in the middle the history generated by the $\mathsf{Eater}$, and below on the right we see the composite history generated by the system as a whole.
\begin{mathpar} 
  \includegraphics[height=3cm,align=c]{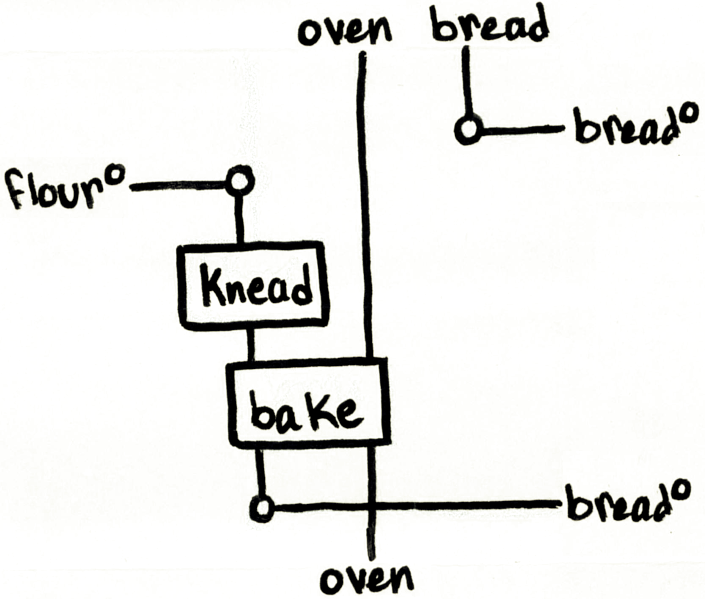}

  \includegraphics[height=3cm,align=c]{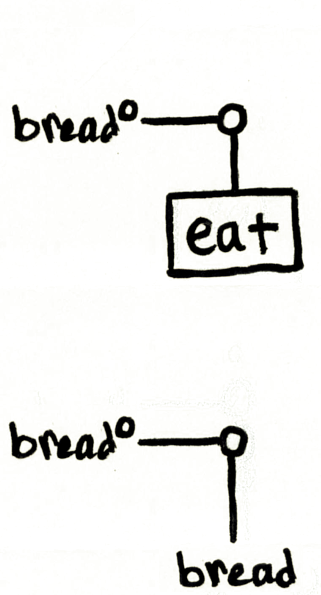}

  \includegraphics[height=3cm,align=c]{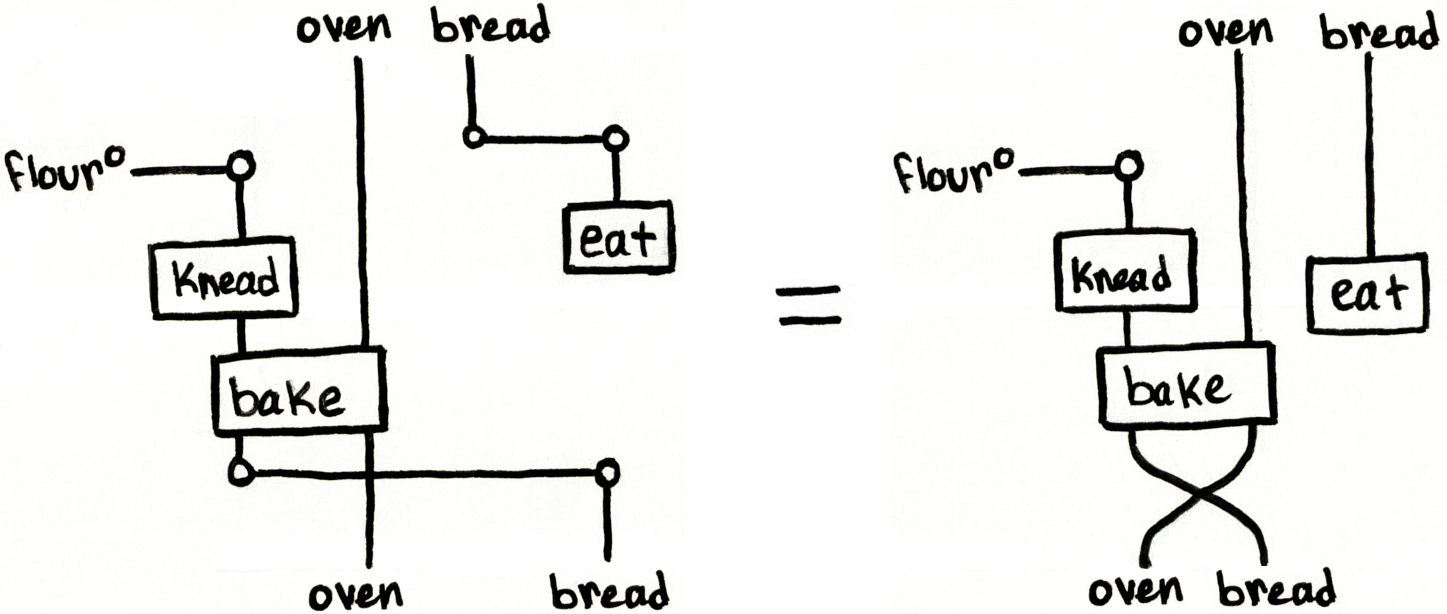}
\end{mathpar}

Situated transition systems are now easily seen to form a category. We record:
\begin{proposition}\label{prop:situated-category}
  Let $\A$ be a resource theory. Then there is a category $\mathsf{S}(\A)$ of situated transition systems, defined as follows:
  \begin{enumerate}[]
  \item \textbf{objects} are $\A$-situated boundaries.
  \item \textbf{arrows} are $\A$-situated transition systems, modulo coherent isomorphism of the underlying spans. That is, for two $\A$-situated transition systems $(R,\phi_R),(S,\phi_S) : (U,\phi_U) \to (V,\phi_V)$, say that $(R,\phi_R) \sim (S,\phi_S)$ in case there exists a reflexive graph isomorphism $\alpha : R \stackrel{\sim}{\to} S$ such that
    \begin{enumerate}[(i)]
    \item $\alpha : R \stackrel{\sim}{\to} S$ is an isomorphism of spans, in the sense that the following diagram commutes:
          \[\begin{tikzcd}
          & R \ar[ld] \ar[dr] \ar[d,"\alpha"] \\
          U & \ar[l] S \ar[r] & V
      \end{tikzcd}\]

    \item $\alpha : R \stackrel{\sim}{\to} S$ preserves material histories, in the sense that there is a natural isomorphism $\iota : \phi_R \stackrel{\sim}{\to} \phi_S \circ \alpha$ (see Remark \ref{remark:natural-enough}).
    \end{enumerate}
      Now an arrow of $\mathsf{S}(\A)$ is a $\sim$-equivalence class of situated transition systems. 
  \item the \textbf{identity} arrow on $(U,\phi_U)$ is given by the identity span $U \stackrel{1_U}{\leftarrow} U \stackrel{1_U}{\rightarrow} U$ and the situation map $\phi_{1_U} : U \to \left< \A \right>$ sends $t : A \to B$ in $U$ to the horizontal identity cell for $\phi_U(t)$. 
  \item \textbf{composition} is as discussed above. 
  \end{enumerate}
\end{proposition}
\qed
\begin{remark}\label{remark:natural-enough}
  In the definition of $\mathsf{S}(\A)$, an equilvalence $(R,\phi_R) \sim (S,\phi_S)$ requires a natural isomorphism $\iota : \phi_R \to \phi_S \circ \alpha$, where $\phi_R$ and $\phi_S \circ \alpha$ are reflexive graph homomorphisms of type $R \to \left< \A \right>$. Natural transformations are defined between \emph{functors}, so the reader would be justified in thinking that we have made a fatal mistake! All is in fact well, as we explain presently. 

  There is a well-known adjunction $F : \mathsf{RGraph} \dashv \mathsf{Cat} : U$ with $F(G)$ being the category of paths in a reflexive graph $G$, and $U(\C)$ being the underlying graph of a category $\C$. Given two reflexive graph homomorphisms $f,g : G \to U(\C)$ define a \emph{natural transformation} $\iota : f \to g$ to consist of a morphism $\iota_A : f(A) \to g(A)$ of $\C$ for each vertex $A$ of $G$ such that for every edge $t : A \to B$ of $G$, $\iota_B \circ f(t) = g(t) \circ \iota_A$ in $\C$. Thus, the definition of natural transformation applies unchanged to reflexive graph homomorphisms whose codomain happens to be a category. Further, applying $F$ to this situation yields a natural transformation in the usual sense. Now $\left< \A \right>$ is clearly the underlying graph of a category, so in particular it makes sense to ask for a natural isomorphism $\iota : \phi_R \to \phi_S \circ \alpha$. Every isomorphism in $\left< \A \right>$ has trivial left and right boundary. We therefore require an isomorphism $\iota_A : \phi_R(A) \stackrel{\sim}{\to} \phi_S(\alpha(A))$ in $\A \cong \bv\,\corner{\A}$ for each vertex $A$ of $R$ such that $\phi_R(t)\iota_B = \iota_A\phi_S(\alpha(t))$ in $\corner{\A}$ for each edge $t : A \to B$ of $R$. 

  Intuitively, isomorphic objects of $\A$ denote the same collection of resources, only orgainzed differently. Understood this way, our notion of equivalence $(R,\phi_R) \sim (S,\phi_S)$ identifies situated transition systems that differ only in the internal organization of their resources. More concretely, asking for strict equality $\phi_R = \phi_S \circ \alpha$ does not result in a monoidal category. We would like $\mathsf{S}(\A)$ to be monoidal, and our notion of equality is just flexible enough to make this the case. 
\end{remark}

\begin{proposition}\label{prop:situated-monoidal}
  If $\A$ is a resource theory then $\mathsf{S}(\A)$ is a monoidal category. 
\end{proposition}
\qed
\section{Compact Closure and Accounting}

In this section we consider the case in which our resource theory $\A$ is compact closed. From the perspective of accountancy, string diagrams over a resource theory are like \emph{ledgers}, recording the material history of the resources they concern \cite{Nes20}. In the partita-doppia (double-entry) method of accounting every change to a ledger must consist of a matching credit (positive change) and debit (negative change), so that the ledger remains \emph{balanced}. This serves as a kind of integrity check: given a ledger we may attempt to \emph{balance} it by matching credits with debits and cancelling them out, and the ledger is well-formed in case all entries may be cancelled in this way.

While the credits and debits of partita-doppia accounting are usually positive and negative integers, the technique applies in the context of any compact closed resource theory. The units $\eta_A : I \to A \otimes A^*$ create matching credits and debits, and the cancellative process of balancing is performed via the counits $\varepsilon_A : A^* \otimes A \to I$. The traditional setting \cite{Ell85} is captured by the compact closed category $\mathbb{Z}$ whose objects are the group of differences construction of the integers and in which there is a morphism between two objects if and only if the corresponding integers are equal \cite{Joy96,Abr05}.

The cells of $\corner{\A}$ with $I$ as their top and bottom boundary are called \emph{horizontal cells}. The horizontal cells of $\corner{\A}$ form a monoidal category $\bh\,\corner{\A}$, with composition given by horizontal composition in $\corner{\A}$ and the tensor product given by vertical composition in $\corner{\A}$. Think of $\bh\,\corner{\A}$ as a category of \emph{exchanges} --- a point of view is developed in \cite{Nes21}. Isomorphic objects of $\bh\,\corner{\A}$ correspond to equivalent exchanges (\cite{Nes21}, Lemma 3). If $\A$ is compact closed we encounter a formal version of the fact that if $\alice$ gives $\bob$ negative five dollars, this is equivalent to $\bob$ giving $\alice$ positive five dollars. More generally, that to get rid of a debit is in many ways the same thing as receiving a credit, and vice-versa. 
\begin{lemma}\label{lem:reversal}
  If $\A$ is compact closed then $A^\circ \cong (A^*)^\bullet$ and $A^\bullet \cong (A^*)^\circ$ in $\bh\,\corner{\A}$.
\end{lemma}
\qed

There is a kind of causal structure present in $\bh\,\corner{\A}$. The corners allow us to bend wires down, but not up, a formal reflection of the fact that I cannot give something away unless I have it. In particular this means that $\mathsf{S}(\A)$ need not be symmetric monoidal: For any $A,B$ there is always a morphism of type $A^\circ \otimes B^\bullet \to B^\bullet \otimes A^\circ$, pictured below on the left, but this is not always an isomorphism. 
\begin{mathpar}
  \includegraphics[height=1cm,align=c]{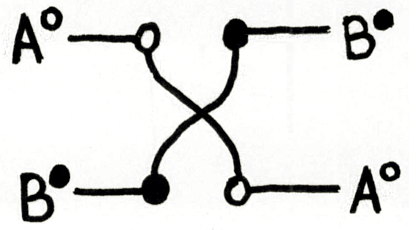}

  \includegraphics[height=1cm,align=c]{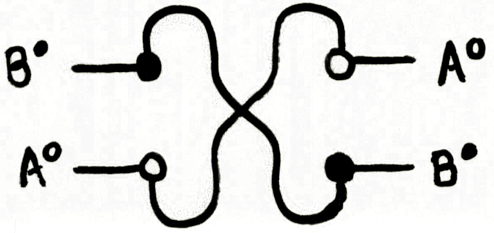}
\end{mathpar}
If our resource theory $\A$ is compact closed, then $\bh\,\corner{\A}$ is symmetric monoidal, with the inverse to the problematic morphism given above on the right. This is a formal reflection of the way that debits allow us to violate causality in everyday life: by incurring a debit I may give something away before I have it. For similar reasons, $\bh\,\corner{\A}$ need not be rigid, but if $\A$ is compact closed then it is. 
\begin{lemma}\label{lem:horiz-compact-closed}
  If $\A$ is compact closed then so is $\bh\,\corner{\A}$.
\end{lemma}
\qed

In fact, if $\A$ is compact closed, then $\mathsf{S}(\A)$ is as well. While we might expect $\mathsf{S}(\A)$ to be compact closed for every $\A$ --- inheriting the compact closed structure of $\spangraph$ --- the geometry of $\bh\,\corner{\A}$ prevents this. Both $\spangraph$ and $\bh\,\corner{\A}$ occur as subcategories of $\mathsf{S}(\A)$, and it seems that structure must be present in both of them in order to manifest in $\mathsf{S}(\A)$. It is interesting that for compact closed resource theories the more flexible compact closed geometry is also present in the category of situated transition systems. Perhaps the use of partita-doppia style debits and credits allows more flexible ``wiring'' of real-world accounting systems than would otherwise be the case. 
\begin{lemma}\label{lem:situated-compact-closed}
  If $\A$ is compact closed, so is $\mathsf{S}(\A)$.
\end{lemma}
\qed

Now, the category $\mathsf{S}(\mathbb{Z})$ of $\mathbb{Z}$-situated transition systems describes systems of partita-doppia accounts in the sense of \cite{Kat98}. The situation maps each state to an integer-valued account balance, and similarly each transition corresponds to a cell of $\corner{\mathbb{Z}}$ with top and bottom boundary the balance of the source and target states, respectively. This ensures that any change in the account balance is reflected by value entering or leaving the system along the boundaries, and vice-versa. Since $\mathbb{Z}$ is compact closed, we obtain an analogue of the main theorem of \cite{Kat98} as a special case of Lemma \ref{lem:situated-compact-closed}, as promised:
\begin{corollary}\label{cor:partita-doppia-corollary}
  $\mathsf{S}(\mathbb{Z})$ is compact closed. 
\end{corollary}
\qed

\section{Conclusions and Future Work}
We have introduced the idea of situating a transition system with boundary in a resource theory and constructed a monoidal category $\mathsf{S}(\A)$ of such systems over an arbitrary resource theory $\A$. Further, we have shown that when $\A$ is compact closed, $\mathsf{S}(\A)$ is also compact closed, generalizing existing work concering systems of partita-doppia accounts \cite{Kat98}. We feel that this in a promising new direction in the study of concurrent systems, and have many ideas for future work. 

If $\A$ is a model of a functional programming language, then an object of $\mathsf{S}(\A)$ can be understood as a very general sort of behavioural type. There is an extensive literature on behavioural types, and we speculate that situated transition systems would be a good way to place this work in the wider context of entire systems. If $\A$ is a model of a ledger system in the sense of \cite{Nes20}, then the material history generated by an $\A$-situated transition system can be seen as a sequence of ledger transactions. It seems that this is relevant to the study of smart contracts, since the ability to transact on the blockchain as they execute is one of their defining features. More ambitiously, we wish to construct compositional models of the systems one encounters in molecular biology, and we imagine that situated transition systems over a resource theory of biomolecules would be a good setting for this. 

It is currently rather painful to specify a situated transition system, and it would be worthwhile to investigate various kinds of syntax that can be given semantics in $\mathsf{S}(\A)$. A promising approach is interpret arrows of $\bh\,\corner{\A}$ as a sort of resource transducer using ideas developed in \cite{Bon19} --- we hope to elaborate on this in a future paper. Finally, ``spancospans'' of reflexive graphs allow us to talk about transition systems with boundary in which the shape of the boundary may change over time \cite{Kat00}. It should be possible to formulate situated transition systems with this capabilty, presumably by working with the intercategory of spancospans \cite{Gra17}.

\bibliographystyle{eptcs}
\bibliography{citations}

\end{document}